\def\Ann{\mathop{\rm Ann}\nolimits}

\def\Ndim{\mathop{\rm N.dim}\nolimits}
\def\width{\mathop{\rm width}\nolimits}

\def\Hom{\mathop{\rm Hom}\nolimits}
\def\id{\mathop{\rm id}\nolimits}
\def\Ext{\mathop{\rm Ext}\nolimits}

\def\LCMo{\mathop{\rm H}\nolimits}

\def\Naturalsign{{\rm l\kern-.23em N}}
\parindent=0pt
\font \normal=cmr10 scaled \magstep0 \font \mittel=cmr10 scaled \magstep1 \font \gross=cmr10 scaled \magstep5
\input amssym.def
\input amssym.tex
\gross
\noindent
Local Homology, Cohen-Macaulayness and Cohen-Macaulayfications
\normal
\bigskip
\bigskip
{\tt Abstract: let }$(R,\goth m)${\tt \ be a local, complete ring, }$X${\tt \ an
artinian }$R${\tt-module of Noetherian \par dimension }d{\tt;let }$x_1,\dots ,x_d\in
\goth m${\tt \ be such that }$0:_X(x_1,\dots ,x_d)R${\tt \ has finite length. Then
}$\LCMo ^{\underline x}_d(X)$ {\tt \ is a finite }$R${\tt -module, providing a positive
answer to a question posed by Tang. As a \par first application of this result
corollary 1 contains a necessary condition for a finite \par module to be
Cohen-Macaulay; secondly we propose a notion of Cohen-Macaulayfication and
\par prove its uniqueness (theorem 3); finally we show that this new notion of
Cohen-\par Macaulayfication is a direct generalization of a notion of
Cohen-Macaulayfication \par introduced by Goto (theorem 4).}
\bigskip
\bigskip
For the sake of completeness we repeat the relevant notions: Throughout this
paper $(R,\goth m)$ is a noetherian, local, complete ring. If $M$ is an
$R$-module and $I$ an ideal of $R$, we denote the $i$-th
local cohomology module of $M$ with support in $I$ by $\LCMo ^i_I(M)$. It is
well-known that $\LCMo ^{\dim (M)}_I(M)$ is artinian for any proper ideal $I$
of $R$ provided $M$ is finitely generated as $R$-module (cf. Melkersson [13]).
\par
There is a theory of local homology modules (cf. Tang [18] and [19]): If $X$ is an
artinian $R$-module and $\underline x=x_1,\dots ,x_r$ is a sequence of
elements in $\goth m$, the $i$-th local homology module $\LCMo ^{\underline
x}_i(X)$ of $X$ with respect to
$\underline x$ is defined by
$$\vtop{\baselineskip=1pt \lineskiplimit=0pt \lineskip=1pt\hbox{lim}
\hbox{$\longleftarrow $} \hbox{$^{^{n\in \bf N}}$}} H_i(K_\bullet (x_1^n,\dots
,x_r^n;X))\ \ ,$$
where $K_\bullet (x_1^n,\dots ,x_r^n;X)$ is the Koszul complex of $X$ with
respect to $x_1^n,\dots ,x_r^n$; then $\LCMo ^{\underline x}_i(\ )$ is an
$R$-linear, covariant functor from artinian $R$-modules to $R$-modules.
\par
We repeat the notions of Noetherian dimension
$\Ndim (X)$ and width of $X$, $\width (X)$: For $X=0$ one puts $\Ndim
(X)=-1$, for $X\neq 0$ $\Ndim (X)$ denotes the least integer $r$ such that
$0:_X(x_1,\dots ,x_r)R$ has finite length for some $x_1,\dots ,x_r\in \goth
m$. Now let $x_1,\dots ,x_n\in \goth m$. $x_1,\dots ,x_n$ is an $X$-coregular
sequence if
$$0:_X(x_1,\dots ,x_{i-1})R\buildrel x_i\over \to 0:X(x_1,\dots ,x_{i-1})R$$
is surjective for $i=1,\dots ,n$. $\width (X)$ is defined as the length
of a (in fact any) maximal $X$-coregular sequence in $\goth m$. Details on
$\Ndim (X)$ and $\width (X)$ can be found in Ooishi [14] and Roberts [15],
here we
cite one general fact: For any artinian $R$-module $X$
$$\width (X)\leq \Ndim (X)<\infty $$
holds and $X$ is co-Cohen-Macaulay if and only if $\width (X)= \Ndim (X)$
holds. Tang has shown ([18], Propostion 2.6) that $\LCMo ^{\dim (M)}_\goth m(M)$ is co-Cohen-Macaulay if $M$
is Cohen-Macaulay and ([18], Theorem 3.1) that $\LCMo ^{x_1,\dots ,x_d}_{\dim (M)}(\LCMo
^{\dim (M)}_\goth m(M))=M$ holds (here $x_1,\dots ,x_d$ is a s. o. p. of
$M$). Tang asks ([18], Remark 3.5)
if $\LCMo ^{\underline x}_d(X)$ is finitely generated if $X$ is
an Artinian $R$-module of N.dimension $d$ and $\underline x=x_1,\dots ,x_d$ is
such that $0:_X\underline x$ has finite length.
\par
We give a positive answer to this question (theorem 1) and draw some consequences
establishing various duality results (theorem 2). As an application we present a necessary
condition for a given finite $R$-module $M$ to be Cohen-Macaulay (corollary 2).
\par
For a given $R$-module $M$ we denote by $D(M)$ the Matlis dual of $M$; for
details on Matlis duality see [3], [4], [5] and [11]. $\underline x$ will
always stand for a sequence $x_1,\dots ,x_d$ in $\goth m$.
\bigskip
\parindent=0pt
{\bf Theorem 1}
\par
Let $(R,\goth m)$ be a noetherian, local, complete ring, $X$ an artinian
$R$-module of N.dimension $d$; let $x_1,\dots ,x_d\in \goth m$ be such that
$0:_X(x_1,\dots ,x_d)R$ has finite length. Then $\LCMo ^{\underline x}_d(X)$ is a finitely generated $R$-module.
\par
Proof:
\par
$x_1,\dots ,x_d$ form a system of parameters for $D(X)$, because
$D(X)/(x_1,\dots ,x_d)D(X)=D(0:_X(x_1,\dots ,x_d)R)$ has finite length and
$\dim (D(X))=\Ndim (X)=d$.
Using Matlis-duality we have
$$\eqalign {\LCMo ^{\underline x}_d(X)&=\LCMo ^{\underline x}_d(D(D(X)))\cr &=\vtop{\baselineskip=1pt \lineskiplimit=0pt \lineskip=1pt\hbox{lim}
\hbox{$\longleftarrow $} \hbox{$^{^{n\in \bf N}}$}} H_d(K_\bullet (x_1^n,\dots ,x_d^n;D(D(X))))\cr &=\vtop{\baselineskip=1pt \lineskiplimit=0pt \lineskip=1pt\hbox{lim}
\hbox{$\longleftarrow $} \hbox{$^{^{n\in \bf N}}$}}D(H^d(K^\bullet (x_1^n,\dots ,x_d^n;D(X))))\cr &=D(\vtop{\baselineskip=1pt \lineskiplimit=0pt \lineskip=1pt\hbox{lim}
\hbox{$\longrightarrow $} \hbox{$^{^{n\in \bf N}}$}}H^d(K^\bullet (x_1^n,\dots ,x_d^n;D(X))))\cr &=D(\LCMo ^d_{(x_1,\dots
,x_d)R}(D(X)))\ \ ,\cr }$$
and the last module is finitely generated because $\LCMo ^d_{(x_1,\dots
,x_d)R}(D(X))$ is artinian.
\bigskip
{\bf Corollary 1}
\par
If $X$ is co-Cohen-Macaulay, then $\LCMo ^{x_1,\dots ,x_d}_d(X)$ is a
Cohen-Macaulay module. In particular if $d=\dim (R)$, $\LCMo ^{x_1,\dots
,x_d}_d(X)$ is a maximal Cohen-Macaulay module.
\par
Proof:
\par
The statements follow from theorem 1 and Tang [18], Remark 3.5.
\bigskip
Let $(R,\goth m)$ be a noetherian, local, complete ring. Let $\cal N$
(resp. $\cal A$) denote the set of isomorphism classes of noetherian (resp. of
artinian) $R$-modules. We have maps $F_1$ and $F_2$ from $\cal N$ to $\cal A$
induced by
$$M\buildrel F_1\over \mapsto \hbox{Matlis dual of }M$$
and
$$M\buildrel F_2\over \mapsto \LCMo ^{\dim (M)}_\goth m(M)$$
For $F_2$ it doesnt make any
difference if we take $\LCMo ^{\dim (M)}_{(x_1,\dots ,x_{\dim (M)})R}(M)$ instead of
$\LCMo ^{\dim (M)}_\goth m(M)$ (for
any system of parameters $x_1,\dots ,x_{\dim (M)}$ of $M$). Similarly we have maps $G_1$ and $G_2$
from $\cal A$ to $\cal N$ induced by
$$X\buildrel G_1\over \mapsto \hbox{Matlis-dual of }X$$
and
$$X\buildrel G_2\over \mapsto \LCMo ^{x_1,\dots ,x_{\Ndim
(X)}}_{\Ndim(X)}(X)$$
(here $x_1,\dots ,x_{\Ndim (X)}$ are such that
$0:_X(x_1,\dots ,x_{\Ndim (X)})R$ has finite length). By
Matlis-duality we have
$$F_1\circ G_1=\id _{\cal A}, G_1\circ F_1=\id _{\cal N}\ \ .$$
From the proof of theorem 1 one sees
$$F_1\circ G_2=F_2\circ G_1=:T$$
and hence
$$G_1\circ F_2=G_2\circ F_1=:T^\prime \ \ ,$$
$$G_2=G_1\circ F_2\circ G_1=G_1\circ T, F_2=F_1\circ G_2\circ F_1=F_1\circ
T^\prime \ \ .$$
\par
{\bf Theorem 2}
\par
Let $(R,\goth m$) be a noetherian, local, complete ring. Let $M$ be a noetherian
and $X$ an artinian $R$-module. Then
\par
(i) If $M$ is Cohen-Macaulay, then $F_2(M)$ is co-Cohen-Macaulay.
\par
(ii) If $M$ is Cohen-Macaulay, then $F_1(M)$ is co-Cohen-Macaulay.
\par
(iii) If $X$ is co-Cohen-Macaulay, then $G_2(M)$ is Cohen-Macaulay.
\par
(iv) If $X$ is co-Cohen-Macaulay, then $G_1(M)$ is Cohen-Macaulay.
\par
Proof:
\par
(ii) and (iv) are easily proved using Matlis-duality theory. (i) is proved by
Tang ([18], Propostion 2.6)  and now (iii) follows from $G_2=G_1\circ F_2\circ G_1$.
\bigskip
Let ${\cal N}_0$ (resp. ${\cal A}_0$) denote the set of isomorphism classes of noetherian Cohen-Macaulay modules (resp. of artinian co-Cohen-Macaulay
modules). Then, by theorem 2, $F_1,F_2,G_1,G_2$ induce maps between ${\cal N}_0$ and ${\cal
A}_0$ in an obvious way. Tang's theorems 3.1 and 3.4 from [18] imply $F_2\circ G_2=\id _{{\cal A}_0}$ and $G_2\circ F_2=\id _{{\cal N}_0}$. We deduce $G_1=G_2\circ
F_1\circ G_2$, $F_1=F_2\circ G_1\circ F_2,T^2=\id ,{T^\prime }^2=\id $ on ${\cal
N}_0$ and ${\cal A}_0$.
\bigskip
As an application there is a necessary condition for a finite module to be Cohen-Macaulay:
\par
{\bf Corollary 2}
\par
(i) Let $\omega _R$ be a dualizing module for $R$ (it exists uniquely up to
isomorphism since $R$ is complete). Assume that $M$ is Cohen-Macaulay. Then
$\Ext ^{\dim (R)-\dim (M)}_R(M,\omega _R)$ is Cohen-Macaulay.
\par
(ii) In particular if there exists an ideal $I$ of $R$ such that $I\subseteq
\Ann _R(M)$, $\dim (R/I)=\dim (M)$ and $R/I$ is Gorenstein, Cohen-Macaulayness
of $M$ implies Cohen-Macaulayness of $\Hom _{\overline R}(M,\overline R)$
(here $\overline R:=R/I$). Such an ideal $I$ exists for example if $R$ itself is Gorenstein.
\par
Proof:
\par
The statements follow from local duality and theorem 2.
\bigskip
In a remark following theorem 2 we have seen $G_2\circ F_2=\id_{{\cal N}_0}$
and $F_2\circ G_2=\id _{{\cal A}_0}$. Now we turn our interest to the
question: What can be said about $G_2\circ F_2$ in general, that is, on $\cal N$?
\bigskip
{\bf Definition}
Let $(R,\goth m)$ be a noetherian, local, complete ring and $M$ a noetherian
(i. e. finitely generated) $R$-module. Let $\tilde M$ be a finitely generated  $R$-module containing $M$
as a submodule. We say $\tilde M$ is a Cohen-Macaulayfication of $M$ if the
following three conditions hold:
\par
(i) $\tilde M$ is Cohen-Macaulay.
\par
(ii) $\dim (\tilde M)=\dim (M)$.
\par
(iii) $\LCMo ^{\dim (M)-1}_\goth m(\tilde M/M)=\LCMo ^{\dim (M)}_\goth
m(\tilde M/M)=0$ (this condition holds for example if $\dim (\tilde M/M)\leq
\dim M-2$).
\bigskip
In the sequel we won't always strictly distinguish between a module $M$ and
its isomorphism class, for reasons of simplicity.
\bigskip
{\bf Theorem 3}
\par
Let $(R,\goth m)$ be a noetherian, local, complete ring and $M$ a noetherian
$R$-module. If $M$ has a Cohen-Macaulayfication, it has (up to an $M$-isomorphism) only one Cohen-Macaulayfication, namely
$(G_2\circ F_2)(M)$.
\par
Proof:
\par
Let $\tilde M$ be a Cohen-Macaulayfication of $M$. We consider the short exact
sequence $0\to M\to \tilde M\to \tilde M/M\to 0$ and its long exact cohomology
sequence induced be applying $\Gamma _\goth m$: Because of condition (iii) of
the definition of a Cohen-Macaulayfication we get a canonical isomorphism $\LCMo ^{\dim
(M)}_\goth m(M)=\LCMo ^{\dim (M)}_\goth m(\tilde M)\buildrel \hbox {(ii)}\over
=\LCMo ^{\dim (\tilde M)}_\goth m(\tilde M)$ and therefore $\tilde M=(G_2\circ
F_2)(\tilde M)=(G_2\circ F_2)(M)$.
\bigskip
{\bf Remark}
\par
Goto (cf. [6]) has shown: If $(A,\goth m)$ is a noetherian, local, $d$-dimensional ring with total
quotient ring $Q(A)$, the following conditions are equivalent:
\par
(i) There is a Cohen-Macaulay ring $B$ between $A$ and $Q(A)$ such that $B$ is
finitely generated as an $A$-module, $\dim (B_\goth n)=d$ for every maximal
ideal $\goth n$ of $B$ and $\goth m\cdot B\subseteq A$.
\par
(ii) $A$ is a Buchsbaum ring (for details on Buchsbaum rings see [17]) and $\LCMo ^i_\goth m(A)=0$ for $i\neq 1,d$.
\par
In this case, if $d\geq 2$, $B$ is uniquely determined and Goto ([6])
calls it  the Cohen-Macaulayfication of $A$.
\bigskip
{\bf Remark}
\par
Cohen-Macaulayfication in our sense is a generalization of Goto's concept of Cohen-Macaulayfication:
\bigskip
{\bf Theorem 4}
\par
Let $(R,\goth m)$ be a noetherian, local, complete ring, and assume that the
Cohen-Macaulayfication $B$ of $R$ (in the sense of Goto) exists. Then $B$ is
also a Cohen-Macaulayfication in our sense.
\par
Proof:
\par
Because of $\goth m\cdot B\subseteq R$ we have $\goth m\cdot (B/R)=0$, which
implies that $B/R$ is a finite-dimensional $R/\goth m$-vector space. Because
of $d=\dim (R)\geq 2$ we must have $\LCMo ^{d-1}_\goth m(B/R)=\LCMo ^d_\goth
m(B/R)=0$.
\bigskip
{\bf Remark}
\par
In particular if $(R,\goth m)$ is a noetherian, local, complete Buchsbaum-ring of
dimension $d\geq 2$ such that $\LCMo ^i_\goth m(R)=0$ for $i\neq 1,d$, the
$R$-module $R$ has a Cohen-Macaulayfication.
\bigskip
{\bf Example}
\par
An easy example is given by $R=k[[x_1,x_2,x_3,x_4]]/(x_1,x_2)\cap
(x_3,x_4)$. In the sense of Goto as well as in our sense $R$ has a
Cohen-Macaulayfication given by $(k[[x_1,\dots ,x_4]]/(x_1,x_2))\oplus
(k[[x_1,\dots ,x_4]]/(x_3,x_4))$; this can be seen either directly or by
remarking that $R$ is a 2-dimensional Buchsbaum ring with $\LCMo ^i_\goth
m(R)=0$ for $i\neq 1,2$.
\bigskip
{\bf Open question}
\par
Let $M$ be a given noetherian $R$-module. Under what conditions is
$G_2(F_2(M))$ a Cohen-Macaulayfication of $M$? Or, equivalently (by theorem 3), when does $M$
have a Cohen-Macaulayfication?
\def\litem{\par\noindent \hangindent=\parindent\ltextindent}
\def\ltextindent#1{\hbox to \hangindent{#1\hss}\ignorespaces}
\bigskip
\mittel
{\bf References}
\normal
\smallskip
\parindent=0.8cm
\litem{1.} Bass, H. On the ubiquity of Gorenstein rings, {\it Math. Z.} {\bf
82}, (1963) 8-28.
\medskip
\litem{2.} Brodmann, M. and Hellus, M. Cohomological patterns of coherent
sheaves over projective schemes, {\it Journal of Pure and Applied Algebra}
{\bf 172}, (2002) 165-182.
\medskip
\litem{3.} Brodmann, M. P. and Sharp, R. J. Local Cohomology, {\it Cambridge
studies in advanced mathematics} {\bf 60}, (1998).
\medskip
\litem{4.} Bruns, W. and Herzog, J. Cohen-Macaulay Rings, {\it Cambridge
University Press}, (1993).
\medskip
\litem{5.} Eisenbud, D. Commutative Algebra with A View Toward Algebraic
Geometry, {\it Springer Verlag}, (1995).
\medskip
\litem{6.} Goto, S. On the Cohen-Macaulayfication of certain Buchsbaum rings,
{\it Nagoya Math. J.} Vol. {\bf 80}, (1980) 107-116.
\medskip
\litem{7.} Grothendieck, A. Local Cohomology, {\it Lecture Notes in
Mathematics, Springer Verlag}, (1967).
\medskip
\litem{8.} Hellus, M. On the set of associated primes of a local cohomology
module, {\it J. Algebra} {\bf 237}, (2001) 406-419.
\medskip
\litem{9.} Hellus, M. On the associated primes of Matlis duals of top local cohomology
modules, to appear in {\it Communications in Algebra}.
\medskip
\litem{10.} Huneke, C. Problems on Local Cohomology, {\it Res.
Notes Math. } {\bf 2}, (1992).
\medskip
\litem{11.} Matlis, E. Injective modules over Noetherian rings, {\it Pacific J. Math.} {\bf 8}, (1958) 511-528.
\medskip
\litem{12.} Matsumura, H. Commutative ring theory, {\it Cambridge
University Press}, (1986).
\medskip
\litem{13.} Melkersson, L. Some applications of a criterion for artinianness of
a module, {\it J. Pure and Appl. Alg.} {\bf 101}, (1995) 293-303.
\medskip
\litem{14.} Ooishi, A. Matlis duality and width of a module, {\it Hiroshima
Math. J.} {\bf 6}, (1976) 573-587.
\medskip
\litem{15.} Roberts R. N. Krull dimension for Artinian modules over quasi local
commutative rings, {\it Quart. J. Math. (Oxford)(3)} {\bf 26}, (1975) 269-273.
\medskip
\litem{16.} Scheja, G. and Storch, U. Regular Sequences and Resultants, {\it AK
Peters}, (2001).
\medskip
\litem{17.} St\"uckrad, J. and Vogel, W. Buchsbaum Rings and Applications. {\it
VEB Deutscher Verlag der Wissenschaften, Berlin, Germany}.
\medskip
\litem{18.} Tang Z. M. Local Homology and Local Cohomology, {\it Algebra
Colloquium} {\bf 11}:4, (2004) 467-476.
\medskip
\litem{19.} Tang Z. M. Local homology theory for Artinian modules, {\it
Comm. Alg.} {\bf 22}, (1994) 2173-2204.

\end